\documentclass[11pt]{article}
\usepackage{amssymb,amsfonts,amsmath,latexsym,epsf,tikz,url}

\newtheorem{theorem}{Theorem}[section]

\newtheorem{conjecture}[theorem]{Conjecture}

\newcommand{\proof}{\noindent{\bf Proof.\ }}
\newcommand{\qed}{\hfill $\square$\medskip}

\begin{document}

\title{Relationship between the distinguishing index, minimum degree  and  maximum degree of graphs}

\author{
Saeid Alikhani  $^{}$\footnote{Corresponding author}
\and
Samaneh Soltani
}

\date{\today}

\maketitle

\begin{center}
Department of Mathematics, Yazd University, 89195-741, Yazd, Iran\\
{\tt alikhani@yazd.ac.ir, s.soltani1979@gmail.com}
\end{center}


\begin{abstract}
Let $\delta$ and $\Delta$ be the minimum and the maximum degree of the vertices of a simple connected graph $G$, respectively.
     The distinguishing index of a graph $G$, denoted by $D'(G)$, is the
least number of labels in an edge labeling of $G$ not preserved by
any non-trivial automorphism.  Motivated by a conjecture by Pil\'sniak (2017)
that implies that for any $2$-connected graph $D'(G) \leq \lceil  \sqrt{\Delta (G)}\rceil +1$, we  prove that for any graph $G$ with $\delta\geq 2$,  $D'(G) \leq \lceil  \sqrt[\delta]{\Delta }\rceil +1$. Also, we show that the distinguishing index of $k$-regular graphs is at most $2$, for any $k\geq 5$.
\end{abstract}


\medskip
\noindent{\bf AMS Subj.\ Class.}: 05C25, 05C15

\section{Introduction }

Let $G=(V,E)$ be a simple connected graph. We use the standard graph notation. In particular, ${\rm Aut}(G)$ denotes the automorphism group of $G$.   For simple connected graph $G$, and $v \in V$,
the \textit{neighborhood} of a vertex $v$ is the set $N_G(v) = \{u \in V(G) : uv \in   E(G)\}$. The \textit{degree of a vertex} $v$ in a graph $G$, denoted by ${\rm deg}_G(v)$, is the number of edges of $G$ incident with $v$. In particular, ${\rm deg}_G(v)$ is the number of neighbours of $v$ in $G$.  We denote by $\delta$ and $\Delta$ the minimum and maximum degrees of the vertices of $G$, respectively. A graph $G$ is \textit{$k$-regular} if ${\rm deg}_G(v) = k$ for all $v \in V$. The \textit{diameter of a graph}  $G$ is the greatest distance between two vertices of $G$, and denoted by ${\rm diam}(G)$.

The \textit{distinguishing index} $D'(G)$ of a graph $G$ is the least number $d$ such that $G$ has an edge labeling  with $d$ labels that is preserved only by the identity automorphism of $G$.  The distinguishing edge labeling was first defined by  Kalinowski and Pil\'sniak \cite{R. Kalinowski and M. Pilsniak} for graphs (was inspired by the well-known distinguishing number $D(G)$ which was defined for general vertex labelings by Albertson and Collins \cite{Albert}). The distinguishing index of some examples of graphs was exhibited in \cite{R. Kalinowski and M. Pilsniak}. For 
 instance, $D'(P_n)=2$ for every $n\geq 3$, and 
 $ D'(C_n)=3$ for $n =3,4,5$,  $D'(C_n)=2$ for $n \geq 6$.  They showed that if $G$ is a connected graph of order $n\geq 3$,  then  $D'(G) \leq \Delta$, unless $G$ is $C_3, C_4$ or $C_5$.  It follows for connected graphs that $D'(G) \geq \Delta$ if and only if $D'(G) =\Delta + 1$ and $G$ is a cycle of length at most five. The equality $D'(G) =
\Delta$ holds for all paths, for cycles of length at least 6, for $K_4$, $K_{3,3}$ and for
symmetric or bisymmetric trees. Also,  Pil\'sniak showed that $D'(G)< \Delta$ for all
other connected graphs.
Pil\'sniak  put forward the following conjecture.

\begin{conjecture} {\rm\cite{nord}}\label{konj}
	If $G$ is a $2$-connected graph, then $D'(G)\leq 1+\lceil  \sqrt{\Delta (G)}\rceil$. 	
\end{conjecture} 
In \cite{upperbound}, we proved that if   $\delta \geq 2$, then $ D'(G)\leq \lceil  \sqrt{\Delta }\rceil +1$,  which proves the conjecture.

\medskip

   Motivated by Conjecture \ref{konj}, in the next section, we prove that for any connected graph  $G$, if  $\delta \geq 2$, then $ D'(G)\leq \lceil  \sqrt[\delta]{\Delta }\rceil +1$.
 Also, in Section 3, we show that the distinguishing index of $k$-regular graphs is at most $2$, for any $k\geq 5$.

\section{An upper bound for $D'(G)$ in terms of $\delta$ and $\Delta$}

In this section, we shall obtain an upper bound for the distinguishing index of graph $G$ in terms of its maximum and minimum degree.   
 For this purpose,  we need some preliminaries.   The friendship graph $F_n$ $(n\geq 2)$ can be constructed by joining $n$ copies of the cycle graph $C_3$ with a common vertex.  The distinguishing index of $F_n$ can be computed by the following result.
 \begin{theorem}{\rm \cite{soltani2}}\label{distfan}
 	Let $a_n=1+27n+3\sqrt{81n^2+6n}$. 
 	For every $n\geq 2$, $$D'(F_n)=\lceil\frac{1}{3} (a_n)^{\frac{1}{3}}+\frac{1}{3(a_n)^{\frac{1}{3}}}+\frac{1}{3}\rceil.$$
 \end{theorem}

Also we need the following theorem: 
\begin{theorem}{\rm\cite{nord}}\label{updisindg}
	Let $G$ be a connected graph that is neither a symmetric nor
	an asymmetric tree. If the maximum degree of $G$ is at least 3, then $D'(G) \leq  \Delta(G) - 1$ unless $G$ is $K_4$ or $K_{3,3}$.
\end{theorem}

\begin{theorem}\label{gooduppdisind}
  For any connected graph $G$, if  $\delta \geq 2$, then $ D'(G)\leq \lceil  \sqrt[\delta]{\Delta }\rceil +1$.
\end{theorem}
\proof 
 If $\Delta \leq 5$, then the result follows from Theorem \ref{updisindg}. So, we suppose that $\Delta \geq 6$. 
   Let $v$ be a vertex of $G$ with the maximum degree $\Delta$.  By Theorem \ref{distfan},  we can label the pendant friendship graph (a subgraph is pendant if it has only
one vertex in common with the rest of a graph) in common with $G$ at $v$ for which $v$ is the central point of the friendship graph, with at most $ \lceil \sqrt{\Delta}\rceil$ labels from label set $\{0, 1, \ldots ,  \lceil \sqrt{\Delta}\rceil \}$, distinguishingly.   If there exists one pendant triangle in common with $G$ at $v$, then we label the two its incident edges to $v$ with 0 and 1, and another edges of the pendant triangle with label 2. 

Let $N^{(1)}(v) =\{v_1, \ldots , v_{|N^{(1)}(v) |}\}$ be the vertices of $G$ at distance one from $v$, except the vertices of pendant friendship or triangle graph in common with $G$ at $v$. Suppose that  $d:=\lceil \sqrt[\delta]{\Delta^{\delta -1}}\rceil -1$ and we continue our labeling by the following steps:
\medskip

Step 1) Since $|N^{(1)}(v)| \leq \Delta$, so we can label the edges $vv_{id +j}$ with label $i$, for $0 \leq i \leq \lceil \sqrt[\delta]{\Delta}\rceil$ and $1 \leq j \leq d$, and we do not use label 0 any more. With respect to the number of incident edges to $v$ with label 0, we conclude that the vertex $v$ is fixed under each automorphism of $G$ preserving the labeling. Also, since the pendant friendship or triangle graph in common with $G$ at $v$ has been labeled distinguishingly, so the vertices of pendant graph  are fixed under each automorphism of $G$ preserving the labeling. Hence, every automorphism of $G$ preserving the labeling must map the set of vertices of $G$ at distance $i$ from $v$ to itself setwise,  for any $1 \leq i \leq {\rm diam}(G)$. We denote the set of vertices of $G$ at distance $i$ from $v$  for any $2 \leq i \leq {\rm diam}(G)$, by $N^{(i)}(v)$.

If $N^{(i)}(v) = \emptyset$, for any $i \geq 2$, then we suppose that  $E_k(v_{jd+k})$ is the set of unlabeled  edges of $G$ incident to the vertex $v_{jd+k}$. For every $0\leq j \leq  \lceil \sqrt[\delta]{\Delta}\rceil$, we can label the elements of each  $E_k(v_{jd +k})$ with labels $\{1, \ldots ,\lceil \sqrt[\delta]{\Delta}\rceil \}$ such that for every  pair of $(E_k(v_{jd+k}), E_k'(v_{jd+k'}))$, where $k\neq k'$, there exist a label $l$, $1\leq l\leq  \lceil \sqrt[\delta]{\Delta}\rceil$, such that the number of label $l$ used for labeling of elements of  $E_k(v_{jd+k})$ and $ E_k'(v_{jd+k'})$ is distinct. Therefore all elements of $N^{(1)}(v)$ is fixed under each automorphism of $G$ preserving the labeling. 
Thus we suppose that $N^{(i)}(v) \neq  \emptyset$, for some $i \geq 2$. 

Now we partition the vertices $N^{(1)}(v)$ to two sets $M_1^{(1)}$ and $M_2^{(1)}$ as follows:
 \begin{equation*}
 M_1^{(1)} =\{x\in N^{(1)}(v)~:~ N(x) \subseteq N(v)\},~~M_2^{(1)}=\{x\in N^{(1)}(v)~:~ N(x) \nsubseteq N(v)\}.
 \end{equation*}
 Thus the sets $M_1^{(1)}$ and $M_2^{(1)}$ are mapped to $M_1^{(1)}$ and $M_2^{(1)}$, respectively, setwise, under each automorphism of $G$ preserving the labeling. For  $0 \leq i \leq \lceil \sqrt[\delta]{\Delta}\rceil$, we set $L_i =\{v_{id +j}~:~ 1 \leq j \leq d\}$. By this notation, we get that for $0 \leq i \leq \lceil \sqrt[\delta]{\Delta}\rceil$, the set $L_i$ is mapped to $L_i$  under each automorphism of $G$ preserving the labeling, setwise. Let the sets $M_{1i}^{(1)}$ and $M_{2i}^{(1)}$ for $0 \leq i \leq \lceil \sqrt[\delta]{\Delta}\rceil$ are as follows:
\begin{equation*}
M_{1i}^{(1)}= M_1^{(1)} \cap L_i,~~ M_{2i}^{(1)}= M_2^{(1)} \cap L_i.
\end{equation*}

It is clear that the sets $M_{1i}^{(1)}$ and $M_{2i}^{(1)}$  are mapped to $M_{1i}^{(1)}$ and $M_{2i}^{(1)}$, respectively, setwise,  under each automorphism of $G$ preserving the labeling. Since  for any $0 \leq i \leq \lceil \sqrt[\delta]{\Delta}\rceil$, we have $|M_{1i}^{(1)}|\leq d$, so we can label all incident edges to each element of $M_{1i}^{(1)}$ with labels $\{1,2, \ldots , \lceil \sqrt[\delta]{\Delta}\rceil \}$, such that for any two vertices of $M_{1i}^{(1)}$, say $x$ and $y$, there exists a label $k$, $1\leq k \leq \lceil \sqrt[\delta]{\Delta}\rceil$, such that the number of label $k$ for the incident edges to $x$ is different from the number of label $k$ for the incident edges to $y$. Hence, it can be deduce that each vertex of  $M_{1i}^{(1)}$ is fixed  under each automorphism of $G$ preserving the labeling, where $0\leq i \leq  \lceil \sqrt[\delta]{\Delta}\rceil$. Thus every vertices of $M_1^{(1)}$ is fixed  under each automorphism of $G$ preserving the labeling. In sequel, we want to label the edges incident to vertices of $M_2^{(1)}$ such that $M_2^{(1)}$ is fixed  under each automorphism of $G$ preserving the labeling, pointwise. For this purpose, we partition the vertices of $M_{2i}^{(1)}$ to the sets $M_{{2i}_j}^{(1)}$, where $1 \leq j \leq \Delta -1$ as follows:
\begin{equation*}
M_{{2i}_j}^{(1)} =\{x\in M_{2i}^{(1)}~:~ |N(x) \cap  N^{(2)}(v)| = j\}.
\end{equation*}

Since the set $N^{(i)}(v)$, for any $i$, is mapped to itself, it can be concluded that $M_{{2i}_j}^{(1)}$  is mapped to itself under each automorphism of $G$ preserving the labeling, for any $i$ and $j$.  Let $M_{{2i}_j}^{(1)} = \{x_{j1}, x_{j2}, \ldots , x_{js_j}\}$. It is clear that $|M_{{2i}_j}^{(1)}| \leq |M_{2i}^{(1)}| \leq d$. Now we consider the two following cases for every $0 \leq i \leq \lceil \sqrt[\delta]{\Delta}\rceil$:

Case 1) Let  $j < \delta -1$ and $\delta \geq 3$. Since   $|M_{{2i}_j}^{(1)}|\leq d $, so we can label all incident edges to each element of $M_{{2i}_j}^{(1)}$ with labels $\{1,2, \ldots , \lceil \sqrt[\delta]{\Delta}\rceil \}$, such that for any two vertices of $M_{{2i}_j}^{(1)}$, say $x$ and $y$, there exists a label $k$, $1\leq k \leq \lceil \sqrt[\delta]{\Delta}\rceil$, such that the number of label $k$ for the incident edges to $x$ is different from the number of label $k$ for the incident edges to $y$. Hence, it can be deduce that each vertex of  $M_{{2i}_j}^{(1)}$ is fixed  under each automorphism of $G$ preserving the labeling, where $1\leq j < \delta -1$.

Case 2) Let  $j \geq \delta -1$. Let $x_{jk}\in M_{{2i}_j}^{(1)}$, and $N(x_{jk}) \cap N^{(2)}(v) =\{x'_{jk1}, x'_{jk2}, \ldots , x'_{jkj}\}$. We assign to the  $j$-tuple $(x_{jk}x'_{jk1}, \ldots , x_{jk}x'_{jkj})$ of edges, a $j$-tuple of labels such that for every $x_{jk}$ and $x_{jk'}$, $1\leq k,k'\leq s_j$, there exists a label $l$ in their corresponding $j$-tuples of labels with different number of label $l$ in their coordinates. For constructing $| M_{{2i}_j}^{(1)}|$ numbers of such $j$-tuples we need, ${\rm min}\{r:~ {j+r-1 \choose r-1}\geq | M_{{2i}_j}^{(1)}| \}$ distinct labels. Since for any $\delta -1 \leq j \leq \Delta -1$, we have
\begin{equation*}
{\rm min}\left\{r:~ {j+r-1 \choose r-1}\geq | M_{{2i}_j}^{(1)}| \right\} \leq {\rm min}\left\{r:~ {j+r-1 \choose r-1}\geq d\right\} \leq \lceil \sqrt[\delta]{\Delta}\rceil,
\end{equation*}
so we need at most $\lceil \sqrt[\delta]{\Delta}\rceil$ distinct labels from label set $\{1, 2, \ldots , \lceil \sqrt[\delta]{\Delta}\rceil\}$ for constructing such $j$-tuples.   Hence, the vertices of $M_{{2i}_j}^{(1)}$, for any $\delta -1 \leq j \leq \Delta -1$, are fixed under each automorphism of $G$ preserving  the labeling.

Therefore, the vertices of $M_{2i}^{(1)}$ for any $0 \leq i \leq \lceil \sqrt{\Delta}\rceil$, and so the vertices of $M_2^{(1)}$ are fixed under each automorphism of $G$ preserving  the labeling. Now, we can get that all vertices of $N^{(1)}(v)$ are fixed. If there exist unlabeled edges of $G$ with the two endpoints in $N^{(1)}(v)$, then we assign them an arbitrary label, say 1.

Step 2) Now we consider $N^{(2)}(v)$. We partition this set such that the vertices of $N^{(2)}(v)$ with the same neighbours in $M_2^{(1)}$, lie in a set. In other words, we can write $N^{(2)}(v) = \bigcup_i A_i$, such that $A_i$ contains that elements of $N^{(2)}(v)$ having the same neighbours in $M_2^{(1)}$, for any $i$. Since all vertices in  $M_2^{(1)}$  are fixed, so  the set $A_i$ is mapped to $A_i$ setwise, under each automorphism of $G$ preserving the labeling. Let $A_i = \{w_{i1}, \ldots , w_{it_i}\}$, and we have 
\begin{equation*}
N(w_{i1}) \cap M_2^{(1)} = \cdots = N(w_{it_i}) \cap M_2^{(1)} = \{v_{i1}, \ldots , v_{ip_i}\}.
\end{equation*}

We consider the two following cases:

Case 1) If for every $w_{ij}$ and $w_{ij'}$ in $A_i$, where $1\leq j,j' \leq t_i$, there exists a $k$, $1\leq k \leq p_i$, for which the label of edges $w_{ij}v_{ik}$ is different from label of edge $w_{ij'}v_{ik}$, then all vertices of $G$ in $A_i$ are fixed under each automorphism  of $G$ preserving the labeling.

Case 2) If there exist  $w_{ij}$ and $w_{ij'}$ in $A_i$, where $1\leq j,j' \leq t_i$, such that for every  $k$, $1\leq k \leq p_i$,  the label of edges $w_{ij}v_{ik}$ and $w_{ij'}v_{ik}$ are the same, then we can make a labeling  such that the  vertices  in $A_i$ have the same property   as Case 1, and so are fixed under each automorphism  of $G$ preserving the labeling, by using at least one of the following actions:
\begin{itemize}
\item By permuting the components of the $j$-tuple of labels assigned to the incident edges to $v_{ik}$ with an end point in $N^{(2)}(v)$,
\item By using a new $j$-tuple of labels, with labels $\{1, 2, \ldots , \lceil \sqrt[\delta]{\Delta}\rceil\}$,  for incident edges to $v_{ik}$ with an end point in $N^{(2)}(v)$, such that the vertices in $M_2^{(1)}$ are fixed under each automorphism of $G$ preserving the labeling, 
\item  By labeling the unlabeled edges of $G$ with the two end points in $N^{(2)}(v)$ which are incident  to the vertices  in $A_i$, 
\item By labeling the unlabeled edges of $G$   which are incident  to the vertices  in $A_i$, and another their endpoint is $N^{(3)}(v)$,
\item By labeling the unlabeled edges of $G$ with the two end points in $N^{(3)}(v)$ for which the end points in $N^{(3)}(v)$ are adjacent  to some of  vertices  in $A_i$.
\end{itemize}

Using at least one of above actions, it can be seen that every two vertices $w_{ij}$ and $w_{ij'}$ in $A_i$ have the property as Case (1). Thus we conclude that all vertices in $A_i$, for any $i$, and so all vertices in $N^{(2)}(v)$, are fixed under each automorphism  of $G$ preserving the labeling.  If there exist unlabeled edges of $G$ with the two endpoints in $N^{(2)}(v)$, then we assign them an arbitrary label, say 1.

By following this method,  in the next step we  partition $N^{(3)}(v)$ exactly by the same method as partition of  $N^{(2)}(v)$ to the sets $A_i$s  in  Step 2, we can make a labeling such that $N^{(i)}(v)$ is fixed pointwise, under each automorphism  of $G$ preserving the labeling, for any $3 \leq i \leq {\rm diam}(G)$.\qed

 By the result obtained  by Fisher and Isaak \cite{fish} and independently by Imrich, Jerebic and Klav\v zar \cite{W Imrich} the distinguishing index of complete bipartite graphs is as  follows. By using Theorem \ref{indcombipar}, we can see that the upper bound of Theorem \ref{gooduppdisind} is sharp for some complete bipartite graphs. 
\begin{theorem}{\rm \cite{fish, W Imrich}}\label{indcombipar}
	Let $p, q, r$ be integers such that $r \geq 2$ and $(r-1)^p <
	q \leq r^p$ . Then
	\begin{equation*}
	D'(K_{p,q}) =\left\{
	\begin{array}{ll}
	r & \text{if}~~ q \leq r^p - \lceil {\rm log}_r p\rceil - 1,\\
	r + 1 & \text{if}~~ q \geq r^p - \lceil {\rm log}_r p\rceil + 1.
	\end{array}\right.
	\end{equation*}
	If $q = r^p -\lceil {\rm log}_r p\rceil$ then the distinguishing index $D'(K_{p,q})$ is either $r$ or $r+1$ and can be computed recursively in $O({\rm log} ∗(q))$ time.
\end{theorem}

 \section{Distinguishing index of regular graphs}
  
By Theorem \ref{gooduppdisind}, we can conclude that the distinguishing index of a $k$-regular graph is at most $3$. In the following we improve this upper bound to $2$.  A \textit{palette}  of a vertex is the set of labels of edges incident to it.
We need the following result to obtain the main result of this section.

\begin{theorem}{\rm\cite{nord}}\label{hamiltopath}
	If $G$ is a  graph of order $n \geq 7$ such that $G$ has a Hamiltonian path, then $D'(G) \leq 2$.
\end{theorem}

\begin{theorem}
Let $G$ be a connected $k$-regular graph of order $n$ with $k \geq 5$. Then $D'(G) \leq 2$.
\end{theorem}
\proof If $k \geq \frac{n-1}{2}$, then it is known that $G$ has a Hamiltonian path, and so $D'(G) \leq 2$, by Theorem \ref{hamiltopath}. Then, we suppose that $5 \leq k < \frac{n-1}{2}$.
 
  Let $v$ be an arbitrary  vertex of $G$, and  $N^{(1)}(v) =\{v_1, \ldots , v_k\}$ be the vertices of $G$ at distance one from $v$.  We state our labeling by the following steps:

Step 1) We label   all incident edges  to $v$ with 1. In our edge labeling of the graph $G$, the vertex $v$ will be the unique vertex with the monochromatic palette $\{1\}$. Hence,   the vertex $v$ is fixed under each automorphism of $G$ preserving the labeling.   Thus, every automorphism of $G$ preserving the labeling must map the set of vertices of $G$ at distance $i$ from $v$ to itself setwise,  for any $1 \leq i \leq {\rm diam}(G)$. We denote the set of vertices of $G$ at distance $i$ from $v$  for any $2 \leq i \leq {\rm diam}(G)$, by $N^{(i)}(v)$. If $N^{(i)}(v) = \emptyset$ for any $i \geq 2$, then $k \geq \frac{n-1}{2}$, which is a contradiction. Thus we suppose that $N^{(i)}(v) \neq \emptyset$ for some $i \geq 2$.

 We can label all incident edges to each element of $N^{(1)}(v) \setminus \{v_1\}$ with labels 1 and 2, such that for any two vertices of $N^{(1)}(v)\setminus \{v_1\}$, say $x$ and $y$, there exists a label $k$, $k = 1,2$, such that the number of label $k$ for the incident edges to $x$ is different from the number of label $k$ for the incident edges to $y$, and also the number of label 2 for the incident edges to each element of  $N^{(1)}(v)\setminus \{v_1\}$ is at least one. Next we label the incident edges to $v_1$ exactly the same as labeling of the incident edges of one of the vertices in $N^{(1)}(v)\setminus \{v_1\}$, say $v_2$. Therefore all vertices in $N^{(1)}(v)$ will also be fixed, except, possibly $v_1$ and $v_2$. To distinguish $v_1$ and $v_2$, we label the incident edges to $v_1$ and $v_2$ which are incident to a vertex in $N^{(2)}(v)$, such that  there exists a label $k = 1,2$, for which  the number of label $k$ for the incident edges to $v_1$ and $v_2$ are distinct. Thus, all vertices in $N^{(1)}(v)$ will be also fixed.

Step 2) Now we consider $N^{(2)}(v)$. We partition this set such that the vertices of $N^{(2)}(v)$ with the same neighbours in $N^{(1)}(v)$, lie in a set. In other words, we can write $N^{(2)}(v) = \bigcup_i A_i$, such that $A_i$ contains that elements of $N^{(2)}(v)$ having the same neighbours in $N^{(1)}(v)$, for any $i$. Since all vertices in  $N^{(1)}(v)$  are fixed, so  the set $A_i$ is mapped to $A_i$ setwise, under each automorphism of $G$ preserving the labeling. Let $A_i = \{w_{i1}, \ldots , w_{it_i}\}$, and we have 
\begin{equation*}
N(w_{i1}) \cap N^{(1)}(v) = \cdots = N(w_{it_i}) \cap N^{(1)}(v) = \{v_{i1}, \ldots , v_{ip_i}\}.
\end{equation*}

We consider the two following cases:

Case 1) If for every $w_{ij}$ and $w_{ij'}$ in $A_i$, where $1\leq j,j' \leq t_i$, there exists a $k$, $1\leq k \leq p_i$, for which the label of edges $w_{ij}v_{ik}$ is different from label of edge $w_{ij'}v_{ik}$, then all vertices of $G$ in $A_i$ are fixed under each automorphism  of $G$ preserving the labeling.

Case 2) If there exist  $w_{ij}$ and $w_{ij'}$ in $A_i$, where $1\leq j,j' \leq t_i$, such that for every  $k$, $1\leq k \leq p_i$,  the label of edges $w_{ij}v_{ik}$ and $w_{ij'}v_{ik}$ are the same, then we can make a labeling  such that the  vertices  in $A_i$ have the same property   as Case 1, and so are fixed under each automorphism  of $G$ preserving the labeling, by using at least one of the following actions:
\begin{itemize}
\item By permuting the  labels assigned to the incident edges to $v_{ik}$ with an end point in $N^{(2)}(v)$,
\item By using a new  labeling  for incident edges to $v_{ik}$ with an end point in $N^{(2)}(v)$, such that the vertices in $N^{(1)}(v)$ are fixed under each automorphism of $G$ preserving the labeling, 
\item  By labeling the unlabeled edges of $G$ with the two end points in $N^{(2)}(v)$ which are incident  to the vertices  in $A_i$, 
\item By labeling the unlabeled edges of $G$   which are incident  to the vertices  in $A_i$, and another their endpoint is $N^{(3)}(v)$,
\item By labeling the unlabeled edges of $G$ with the two end points in $N^{(3)}(v)$ for which the end points in $N^{(3)}(v)$ are adjacent  to some of  vertices  in $A_i$.
\end{itemize}

Using at least one of above actions, it can be  concluded that all vertices in $A_i$, for any $i$, and so all vertices in $N^{(2)}(v)$, are fixed under each automorphism  of $G$ preserving the labeling.  If there exist unlabeled edges of $G$ with the two endpoints in $N^{(2)}(v)$, then we assign them an arbitrary label, say 2.

By following this method,  in the next step we  partition $N^{(3)}(v)$ exactly by the same method as partition of  $N^{(2)}(v)$ to the sets $A_i$s  in  Step 2, we can make a labeling such that $N^{(i)}(v)$ is fixed pointwise, under each automorphism  of $G$ preserving the labeling, for any $3 \leq i \leq {\rm diam}(G)$.\qed

\end{document}